\documentclass{amsart}

\usepackage{url}
\usepackage{graphicx}
\newtheorem{thm}{Theorem}[section]
\newtheorem{cor}{Corollary}[section]

\theoremstyle{definition}

\theoremstyle{remark}
\newtheorem{rem}[thm]{Remark}

\begin{document}
%\linespread{1.25}
\title[On the stability of the quaternion projective space]{On the stability of the quaternion projective space}

\author[C.-D. Neac\c{s}u]{Crina-Daniela Neac\c{s}u}

\date{}
\maketitle

\abstract
The aim of this note is to prove that the index of the identity map on a quaternion space form of constant quaternion sectional curvature $c$ is zero, provided that $c<0$. As an immediate consequence, it is established that any compact quaternion space form of negative quaternion sectional curvature is stable and it is emphasized that, on the contrary (but in agreement with some known results), any quaternion projective space is unstable.
\\ \\
{\bf Keywords:} Harmonic map; stability; quaternion space form; identity map.\\ \\
{\bf Mathematics Subject Classification (MSC2020):} 58E20.
\endabstract

\section{Introduction}

The concept of harmonic map between two Riemannian manifolds $(M,g)$ and $(N,h)$ was really introduced in its full generality
in 1964 by Eells and Sampson \cite{ES} as a generalisation of the notion of geodesic,  %the authors
showing that in some geometric settings an arbitrary function can be deformed into an harmonic map. It is known that these maps, that are nothing but critical points of the Dirichlet energy, are of interests in theoretical physics where they appear under the name of generalized sigma models (see \cite{Eich}). Note that first sigma model was introduced in physics by Zakharov and Mikhailov \cite{ZM} 14 years after the publication of the fundamental work \cite{ES}.
The importance of harmonic maps for some different physical phenomena was first outlined by Misner in \cite{Misner}.
Later, in an excellent survey article on harmonic maps from the viewpoint of physics,
S\'{a}nchez \cite{San} describes, in a unified framework, sigma models, Yang-Mills fields and the Einstein equations. It is also outlined that harmonic maps provide us a very natural setting for discussing the Hamiltonian formulation of field theories, such as general relativity (see also \cite{Sci}). For some recent results on the harmonicity of various types of maps defined from Riemannian manifolds endowed with different geometric structures, see \cite{ASah1,ASah2,BG,Kau,Sahin,Vilcu,Wani}.

Concerning the study of harmonic maps, one of the subjects of high interest common to both mathematicians and physicists consists in investigating the stability property of such maps (see, e.g., \cite{BW,TOT,URK}). Recall that a harmonic map $h:(M,g_M)\rightarrow (N,g_N)$ is called stable if the second variation of the  Dirichlet energy is nonnegative
for any smooth variation of the map. Equivalently, a harmonic map $h:(M,g_M)\rightarrow (N,g_N)$ is stable if the index of $h$, denoted by $\mathfrak{i}(h)$, is zero.

It is obvious that the simplest example of harmonic map is provided by the identity map of a Riemannian space $(M,g_M)$ (see \cite{Chen}). If this map $1_M:(M,g_M)\rightarrow(M,g_M)$ is stable, then $(M,g_M)$ is called stable. Otherwise, $(M,g_M)$ is said to be unstable.
Although the identity map $1_M$ has an extremely simple form, investigating its stability is a very challenging subject and numerous interesting results have been obtained over time. For example, in a series of three works published over 7 years, Rehman investigated the stability of generalized Sasakian space forms, $S$-space forms and  Lorentzian para Sasakian manifolds \cite{Reh15, Reh19,Reh22}. Moreover, Gherghe and V\^{\i}lcu \cite{GV2023} established a criterion for the stability of the identity map on locally conformal almost cosymplectic manifolds of
pointwise constant $\phi$-holomorphic sectional curvature. Very recently, the author of the present paper proved that a compact $T$-space form $M$ having non-positive constant $\phi$-sectional curvature $c$ is unstable if the first eigenvalue of the Laplace-Beltrami operator admits a suitable upper bound in terms of $c$ and dimension of $M$
(see \cite{CriCriGreieras}). For other interesting results concerning the stability of harmonic maps on  manifolds equipped with various geometric structures, reference is made to \cite{Burns,Gherghe13,GIP,IMV08,Ohnita,Smith}.

Motivated by all these results, we are going to investigate the stability of the quaternion space forms $M^n(c)$ of quaternion dimension $n$ and  constant quaternion sectional curvature $c$. The main result of the paper is the following.

\begin{thm}\label{Th1}
  The index of the identity map of a compact quaternion space form $M^n(c)$ is zero, i.e.
  $\mathfrak{i}(1_{M^n(c)})=0$, provided that $c<0$.
\end{thm}

As a direct outcome of the above theorem, one arrives at the following consequence.

\begin{cor}
  Any compact quaternion space form of negative quaternion sectional curvature is stable.
\end{cor}

Note that for complex projective spaces, it was proved by Burns, Burstall, De Bartolomeis and Rawnsley a more general result, namely that any compact K\"{a}hler manifold is stable (see  \cite[Lemma 3]{Burns}). Thus, it follows in particular that the complex projective space $P^n(\mathbb{C})(c)$ of constant holomorphic sectional curvature $c$ is stable.

\section{Preliminaries}

\subsection{Quaternion K\"{a}hler manifolds}

Let $(M,g_M)$ be an $m$-dimensional Riemannian space equipped with subbundle $\Lambda\subset End(TM)$ of rank 3, with local basis $\lbrace{J_i}\rbrace_{i=1,2,3}$ such that:
    \begin{equation}\label{cm}
    g_M(J_i\cdot,J_i\cdot)=g_M(\cdot,\cdot)
    \end{equation}
and
       \begin{equation}\label{qi}
       J_i^2=-{\rm I},\quad
       J_{i}J_{i+1}=-J_{i+1}J_{i}=J_{i+2},
       \end{equation}
for $i=1,2,3$, where ${\rm I}$ stands for the identity tensor field of
type (1, 1) on $M$ and the indices are taken from the set $\{1,2,3\}$ modulo 3.
Then the triple $(M,\Lambda,g_M)$ is called \emph{an almost
quaternion Hermitian manifold}. It is not difficult to prove that such manifolds are of dimension divisible by $4$, since one can construct a local orthonormal frame on $M$ in the form \[\{E_1,...,E_n,J_1E_1,...,J_1E_n,J_2E_1,...,J_2E_n,
J_3E_1,...,J_3E_n\}.
\]
Hence, the dimension of $M$ is indeed $m=4n$, where $n\geq 1$. Note that $m$ is called the \emph{real dimension} of the almost
quaternion Hermitian manifold $M$, while $n=\frac{m}{4}$ is said to be the \emph{quaternion dimension} of $M$.
Recall that if $\Lambda$ is a parallel subbundle with respect to the Levi-Civita
connection $\nabla^{g_M}$ of $g_M$, then the triple
$(M,\Lambda,g_M)$ is called a {\it quaternion
K\"{a}hler manifold.}

Suppose now that $(M,\Lambda,g_M)$  is a quaternion
K\"{a}hler manifold and $X\neq 0$ is a vector field on
$M$, then the 4-plane spanned by $\{X,J_1X,J_2X,J_3X\}$, denoted by $\mathcal{Q}(X)$,
is said to be a \emph{quaternion 4-plane}. If $\Pi\subset \mathcal{Q}(X)$ is a 2-plane, then $\Pi$ is said to be a \emph{quaternion plane}. The sectional curvature of a
quaternion plane is termed as a \emph{quaternion sectional curvature}. A
quaternion K\"{a}hler manifold is called a \emph{quaternion space form} if its
quaternion sectional curvatures are constant.
For a quaternion space form of constant quaternion sectional curvature $c$ and quaternion dimension $n$, usually denoted by $M^n(c)$,
the Riemannian curvature tensor $R$ satisfies \cite{BCU}
%     \begin{equation}\begin{aligned}\label{6}
%       R^M(X,Y)Z= &\frac{c}{4}\lbrace g_M(Z,Y) X-
%      g_M(X,Z) Y+\sum\limits_{i=1}^3
%        [g_M(Z,J_i Y) J_i X-\\
%       &-g_M(Z,J_i X) J_i Y+
%        2g_M(X,J_i Y) J_i Z]\rbrace
%    \end{aligned}  \end{equation}
%where $X,Y,Z$ are vector fields on $M$ and
%$\lbrace{J_1,J_2,J_3}\rbrace$ is a local basis of $\Lambda$.
%Using now \eqref{6}, we obtain
     \begin{eqnarray}\label{6bis}
       g_M(R^M(X,Y)Z,U)&= &-\frac{c}{4}\lbrace g_M(X,Z) g_M(Y,U)-g_M(Z,Y) g_M(X,U)
      \nonumber\\
      &&+\sum\limits_{\alpha=1}^3
        [g_M(X,J_\alpha Z) g_M(Y,J_\alpha U)\nonumber\\
       &&-g_M(U,J_\alpha X) g_M(J_\alpha Y,Z)+
        2g_M(X,J_\alpha Y) g_M(J_\alpha U,Z]\rbrace,
    \end{eqnarray}
for any vector fields $X,Y,Z,U$ on $M$.

It is a well-known result that a quaternion space form is locally congruent to either a quaternion projective space $P^n(\mathbb{H})(c)$ of quaternion sectional curvature $c>0$, a quaternion Euclidean space $\mathbb{H}^n$ of null quaternion sectional curvature or a quaternion hyperbolic space $\mathbb{H}H^n(c)$ of quaternion sectional curvature $c<0$ (see, e.g., \cite{AM,Alodan,Aquib,ORTPER}. %In the following, we will focuss our study on the quaternion projective space $\mathbb{H}P^n(c)$, which is known to be a compact quaternion K\"{a}hler symmetric space (see \cite{BES}).

\subsection{Harmonic maps}

Suppose $h:(M,g_M)\rightarrow (N,g_N)$ is a smooth map between
two Riemannian spaces $(M,g_M)$ and $(N,g_N)$. Then the \emph{second fundamental form} $\alpha_h$ of $h$ is defined as
\[\alpha_h(X,Y)=\widetilde{\nabla}_X h_*Y-h_*\nabla_X Y,\]
for any
vector fields $X,Y$ on $M$, where $\nabla$ is the Levi-Civita
connection of $M$ and $\widetilde{\nabla}$ is the pullback of the
connection $\nabla'$ of $N$ to the induced vector bundle
$h^{-1}(TN)$, i.e. \[\widetilde{\nabla}_Xh_*Y=\nabla'_{h_*X}h_*Y.\] The
tension field $\tau(h)$ of the smooth map $h$ is defined as the trace of
$\alpha_h$, that is
    \begin{equation}\label{5}
    \tau(h)_p=\sum_{i=1}^m\alpha_h(E_i,E_i),
    \end{equation}
where $\{E_1,...,E_m\}$ is a local orthonormal frame of $T_pM$,
$p\in M$.

Moreover, if $M$ is compact, then the \emph{Dirichlet energy} of $h$ is
defined by \[\mathcal{E}(h)=\int_M e(h)\vartheta_{g_M},\] where $\vartheta_{g_M}$ is the
canonical measure associated with the metric $g_M$ and
\[e(h)_p=\frac{1}{2} \mathrm{tr}(h^*g_{N})_p,\quad \forall p\in M.\]

We recall that the map $h$ is said to be harmonic if for any smooth variation $\{h_t\}_{t\in(-\epsilon,\epsilon)}$ of $h$, with $h_0=h$, one has
\[
\frac{d}{dt}\mathcal{E}_t|_{t=0}=0.
\]
Equivalently, $h$ is harmonic
if and only if $\tau(h)_p=0$, for all $p\in M$.

If $\{h_{s,t}\}_{s,t\in(-\epsilon,\epsilon)}$ is a smooth
variation of $h$ with two parameters $s$ and $t$ such that $h_{0,0}=h$, then the \emph{Hessian}
of $h$ is defined as:
\[\mathrm{Hess}_h(V,W)=\frac{\partial^2}{\partial s\partial
t}(\mathcal{E}(h_{s,t}))|_{(s,t)=(0,0)},\] where
$V,W\in\Gamma(h^{-1}(TN))$ are the associated variational vector
fields, that is
\[
V=\frac{\partial}{\partial s}(h_{s,t})|_{(s,t)=(0,0)}
\]
and
\[
W=\frac{\partial}{\partial t}(h_{s,t})|_{(s,t)=(0,0)}.
\]

For a harmonic map $h$, the index of $h$ is a natural number denoted by $\mathfrak{i}(h)$ and defined as
the dimension of the largest subspace of $\Gamma(h^{-1}(TN))$ on
which the Hessian of $h$ is negative definite. A harmonic map is called stable if $\mathfrak{i}(h)=0$. Otherwise, $h$ is called
unstable. A key ingredient in study the stability of harmonic maps is provided by
the second variation formula due to Mazet and Smith \cite{Mazet,Smith}
    \begin{equation}\label{5a}
    \mathrm{Hess}_h(V,W)=\int_M g_N(\mathcal{J}_h(V),W)\vartheta_{g_M},
    \end{equation}
where $\emph{J}_h$ is the Jacobian operator of the map $h$ defined in \cite{BW} by
\begin{equation}\label{eq8}
  \mathcal{J}_hV=-\sum_{i=1}^{m}\left(\widetilde{\nabla}_{E_i}\widetilde{\nabla}_{E_i}-\widetilde{\nabla}_{\nabla_{E_i}E_i}\right)V-
   \sum_{i=1}^{m}R^N(V,h_{*}E_i)h_{*}E_i,
\end{equation}
for any $V\in\Gamma(h^{-1}(TN))$, where $\{E_1,...,E_m\}$ is a local orthonormal frame on $M$ and $R^N$ denotes the Riemannian curvature tensor of $N$.

Using the following well-known formula of the rough Laplacian $\bar{\Delta}_h$ of $h$:
 \begin{equation}\label{eq9}
 \bar{\Delta}_hV=-\sum_{i=1}^{m}\left(\widetilde{\nabla}_{E_i}\widetilde{\nabla}_{E_i}-\widetilde{\nabla}_{\nabla_{E_i}E_i}\right)V,
 \end{equation}
one derives from \eqref{eq8} the next expression for $J_h$:
\begin{equation}\label{eq10}
  \mathcal{J}_hV=\bar{\Delta}_hV-
   \sum_{i=1}^{m}R^N(V,h_{*}E_i)h_{*}E_i.
  \end{equation}

\section{Proof of Theorem \ref{Th1}}

Suppose $M=M^n(c)$ is a quaternion  space form of quaternion sectional curvature $c$ and quaternion dimension $n$, equipped with the compatible  Riemannian metric $g_M$. We consider the identity function of this space $1_M:M\rightarrow M$, defined by $1_M(p)=p$, for any  $p\in M$. Applying the second variation formula \eqref{5a} for $1_M$  and considering \eqref{eq8}, we deduce:
\begin{eqnarray}\label{He1}
    \mathrm{Hess}_{1_M}(V,V)&=&\int_Mg_M(\mathcal{J}_{1_M}V,V)\vartheta_{g_M}\nonumber\\
 &=&-\int_M\sum_{i=1}^{n}g_M((\widetilde{\nabla}_{E_i}\widetilde{\nabla}_{E_i}-\widetilde{\nabla}_{\nabla_{E_i}E_i})V,V)\vartheta_{g_M}\nonumber\\
&&-\int_M\sum_{i=1}^{n}g_M(R^M(V,E_i)E_i,V)\vartheta_{g_M}\nonumber\\
&&-\int_M\sum_{\alpha=1}^{3}\sum_{i=1}^{n}g_M((\widetilde{\nabla}_{J_{\alpha}E_i}\widetilde{\nabla}_{J_{\alpha}E_i}
-\widetilde{\nabla}_{\nabla_{J_{\alpha}E_i}{J_{\alpha}E_i}})V,V)\vartheta_{g_M}\nonumber\\
&&-\int_M\sum_{\alpha=1}^{3}\sum_{i=1}^{n}g_M(R^M(V,{J_{\alpha}}E_i){J_{\alpha}}E_i,V)\vartheta_{g_M},
 \end{eqnarray}
where
$\{E_1,...,E_n,J_1E_1,..,J_1 E_n,J_2E_1,..,J_2 E_n,J_3E_1,..,J_3 E_n\}$  is a local orthonormal frame on $M$.

But the expression of the rough Laplacian  $\bar{\Delta}_h(V)$ with respect to the above frame is
 \begin{eqnarray}\label{eq9bis}
 \bar{\Delta}_hV&=&-\sum_{i=1}^{n}(\widetilde{\nabla}_{E_i}\widetilde{\nabla}_{E_i}-\widetilde{\nabla}_{\nabla_{E_i}E_i})V\nonumber\\
 &&-\sum_{\alpha=1}^{3}\sum_{i=1}^{n}(\widetilde{\nabla}_{J_{\alpha}E_i}\widetilde{\nabla}_{J_{\alpha}E_i}
-\widetilde{\nabla}_{\nabla_{J_{\alpha}E_i}{J_{\alpha}E_i}})V.
 \end{eqnarray}

Using \eqref{eq9bis} in \eqref{He1} we obtain immediately that the Hessian of $1_M$ can be rewritten as
\begin{eqnarray}\label{eq101}
    \mathrm{Hess}_{1_M}(V,V)&=&\int_Mg_M(\bar{\Delta}_{1_M}V,V)\vartheta_{g_M}\nonumber\\
    &&+\int_M\sum_{i=1}^{n}g_M(R^M(E_i,V)E_i,V)\vartheta_{g_M}\nonumber\\
    &&+\int_M\sum_{\alpha=1}^{3}\sum_{i=1}^{n}g_M(R^M(J_{\alpha}E_i,V)J_{\alpha}E_i,V))\vartheta_{g_M},
    \end{eqnarray}
where we used the antisymmetry of the Riemannian curvature tensor $R^M$ in relation to the first two arguments.

Replacing now $X=E_i$, $Y=V$, $Z=E_i$ and $U=V$ in \eqref{6bis},
in view of the fact that
\[
g(E_i,E_i)=1,\quad
g(E_i,J_1 E_i)=g(E_i,J_2 E_i)=g(E_i,J_3 E_i)=0,
\]
and taking into account of \eqref{cm} and \eqref{qi},
we find
\begin{eqnarray}\label{cri1}
 g_M(R^M(E_i, V) E_i, V)&=& -\frac{c}{4}\{g_M(V,V)-g_M(V, E_i)g_M(E_i,V)\nonumber\\
  && +3\sum_{\beta=1}^{3}[g_M(V, J_{\beta}E_i)g_M(V,J_{\beta}E_i)]\}.
\end{eqnarray}

Similarly, replacing $X=J_{\alpha}E_i$, $Y=V$, $Z=J_{\alpha}E_i$ and $U=V$ in \eqref{6bis}, we get using the same arguments in computation that
\begin{eqnarray}\label{cri2}
 g_M(R^M(J_{\alpha} E_i,V)J_{\alpha} E_i, V)&=&-\frac{c}{4}\{g_M(V,V)-g_M(V, J_{\alpha} E_i)g_M(J_{\alpha} E_i,V)\nonumber\\
  && +3\sum_{\beta=1}^{3}[g_M(V, J_{\beta}J_{\alpha}E_i)g_M(V,J_{\beta}J_{\alpha}E_i)]\}.
\end{eqnarray}

Next, a use of \eqref{cri1} and \eqref{cri2} in \eqref{eq101} leads to

\begin{eqnarray}
 \mathrm{Hess}_{1_M}(V,V)&=&\int_Mg_M(\bar{\Delta}_{1_M}V,V)\vartheta_{g_M}\nonumber\\
&&-\int_M\sum_{i=1}^{n} \frac{c}{4}\{g_M(V,V)-g_M(V,  E_i)g_M(E_i,V)\nonumber\\
  && +3\sum_{\beta=1}^{3}[g_M(V, J_{\beta}E_i)g_M(V,J_{\beta}E_i)]\}\vartheta_{g_M}\nonumber\\
&&-\int_M\sum_{\alpha=1}^{n}\sum_{i=1}^{n} \frac{c}{4}\{g_M(V,V)-g_M(V, J_{\alpha} E_i)g_M(J_{\alpha}E_i,V)\nonumber\\
&&+3\sum_{\beta=1}^{3}[g_M(V, J_{\beta}J_{\alpha}E_i)g_M(V,J_{\beta}J_{\alpha}E_i)]\}\vartheta_{g_M}.\nonumber
\end{eqnarray}

Adding and combining the terms of the same type in the above equation, we derive that the Hessian of $1_M$ can be expressed as

\begin{eqnarray}\label{co}
 \mathrm{Hess}_{1_M}(V,V)&=&\int_Mg_M(\bar{\Delta}_{1_M}V,V)\vartheta_{g_M}-nc\int_Mg_M(V,V)\vartheta_{g_M}\nonumber\\
&&+ \frac{c}{4}\int_M[\sum_{i=1}^{n}g_M(V,E_i)g_M(E_i,V)\nonumber\\
&&\ \ \ \ \ \ +\sum_{\alpha=1}^{3}\sum_{i=1}^{n}g_M(V,J_{\alpha}E_i)g_M(J_{\alpha}E_i,V)]\vartheta_{g_M}\nonumber\\
&&- \frac{3c}{4}\int_M[\sum_{\beta=1}^{3}g_M(V,J_{\beta}E_i)g_M(V,J_{\beta}E_i)\nonumber\\
&&\ \ \ \ \ \ +\sum_{\alpha=1}^{n}\sum_{\beta=1}^{3}g_M(V, J_{\beta}J_{\alpha}E_i)g_M(V, J_{\beta}J_{\alpha}E_i)]\vartheta_{g_M}.\nonumber\\
\end{eqnarray}

But using now some standard argument from linear algebra, we can remark easily that

\begin{eqnarray}\label{co1}
&&\sum_{i=1}^{n}g_M(V,E_i)g_M(E_i,V)
+\sum_{\alpha=1}^{3}\sum_{i=1}^{n}g_M(V,J_{\alpha}E_i)g_M(J_{\alpha}E_i,V)\nonumber\\
&&=\sum_{i=1}^{n}{g_M(V,g_M(E_i,V)E_i)}+\sum_{\alpha=1}^{3}g_M(V,g_M(J_{\alpha}E_i, V)J_{\alpha}E_i)\nonumber\\
&&=g_M(V,\sum_{i=1}^{n}g_M(E_i,V)E_i+\sum_{\alpha=1}^{3}\sum_{i=1}^{n}g_M(J_{\alpha}E_i, V)J_{\alpha}E_i)\nonumber\\
&&=g_M(V, V),
\end{eqnarray}
and in a similar way we derive
\begin{eqnarray}\label{co2}
&&\sum_{\beta=1}^{3}g_M(V,J_{\beta}E_i)g_M(V,J_{\beta}E_i)+\sum_{\alpha=1}^{n}\sum_{\beta=1}^{3}g_M(V, J_{\beta}J_{\alpha}E_i)g_M(V, J_{\beta}J_{\alpha}E_i)\nonumber\\
&&=\sum_{\beta=1}^{3}g_M(V,g_M(V, J_{\beta}E_i)J_{\beta}E_i)+\sum_{\alpha=1}^{n}g_M(V,g_M(V, J_{\beta}J_{\alpha}E_i)J_{\beta}J_{\alpha}E_i)]\nonumber\\
&&=\sum_{\beta=1}^{3}g_M(J_{\beta}V,g_M(J_{\beta}V,E_i)E_i)+\sum_{\alpha=1}^{n}g_M( J_{\beta}V, g_M(J_{\beta}V, J_{\alpha}E_i)J_{\alpha}E_i)]\nonumber\\
&&=\sum_{\beta=1}^{3}g_M(J_{\beta}V,J_{\beta}V)\nonumber\\
&&=\sum_{\beta=1}^{3}g_M(V, V)\nonumber\\&&=3g_M(V,V).
\end{eqnarray}

So, replacing \eqref{co1} and \eqref{co2} in \eqref{co}, we obtain immediately that
\begin{eqnarray}\label{coconut}
 \mathrm{Hess}_{1_M}(V,V)&=&\int_Mg_M(\bar{\Delta}_{1_M}V,V)\vartheta_{g_M}\nonumber\\
 &&-(n+2)c\int_Mg_M(V, V)\vartheta_{g_M}.
\end{eqnarray}

But it is known (see, e.g., Theorem 1 in \cite{Reh19}) that
 \begin{eqnarray}\label{Reh}
 \int_Mg_M(\bar{\Delta}_{1_M}V,V)\vartheta_{g_M}= \int_Mg_M(\widetilde{\nabla}V,\widetilde{\nabla}V)\vartheta_{g_M}.
\end{eqnarray}

Consequently, inserting \eqref{Reh} in \eqref{coconut}, we derive
\begin{eqnarray}
 \mathrm{Hess}_{1_M}(V,V)&=&\int_Mg_M(\widetilde{\nabla}V,\widetilde{\nabla}V)\vartheta_{g_M}\nonumber\\
&&-(n+2)c\int_Mg_M(V, V)\vartheta_{g_M},
\end{eqnarray}
and therefore it is clear that
\[
\mathrm{Hess}_{1_M}(V,V)\geq{0},
\]
if the constant quaternion sectional curvature $c$ is negative.

 Thus, we can conclude that the Hessian of the identity map $1_M$ is positive definite, which implies that the index of the identity map of quaternion space form $M^n(c)$ is zero, if $c<0$.
Hence the proof of Theorem \ref{Th1} is now complete and as an immediate consequence we have that any compact quaternion space form $M^n(c)$ is stable, provided that  $c<0$.

\section{Final remarks}

\begin{rem}
  Note that the standard example of a quaternionic space form of negative quaternionic sectional curvature is the quaternionic hyperbolic space $\mathbb{H}H^n$. However, this space is non-compact. Note that compact quaternion space forms with $c<0$ do exist, but only as quotients of the quaternion hyperbolic space $\mathbb{H}H^n$ by discrete cocompact groups of isometries.
An explicit example of a compact quaternion space form of negative quaternion sectional curvature can be constructed as follows. If we take
$\Gamma $ be a discrete, cocompact, torsion-free subgroup (i.e., a cocompact lattice) of the isometry group of $\mathbb{H}H^n$, then the quotient
$M = \mathbb{H}H^n / \Gamma$ is a compact manifold with constant negative quaternion sectional curvature.
\end{rem}

\begin{rem}
Concerning the stability of the quaternion projective space $P^n(\mathbb{H})$
we can remark that this space is unstable. This result follows as a consequence of \cite[Proposition 2.11]{Smith}. According to this Proposition, a closed oriented Einstein manifold $(M,g)$ with Einstein constant $\mathcal{C}$ is stable if and only if the first eigenvalue $\lambda_1$ of the Laplace-Beltrami operator satisfies $\lambda_1\geq 2\mathcal{C}$. But it is known that the quaternion projective space $P^n(\mathbb{H})$ is an Einstein space with Einstein constant $4(n+2)$ (see \cite{BES}), while the first eigenvalue of the Laplace-Beltrami operator is $\lambda_1=8(n+1)$ (see \cite{Jost}). Hence we have
\[
\lambda_1=8(n+1)< 8(n+2)=2\mathcal{C}.
\]
Thus, applying \cite[Proposition 2.11]{Smith}, it follows that $P^n(\mathbb{H})$ is unstable. This result is in agreement with \cite[Theorem 1]{Ohnita}.
\end{rem}

\section*{Data Availability Statement}
Data availability is not applicable to this article as no new data were created or analysed in this study.

\section*{Conflicts of Interest}
The author declares no conflict of interest.

\section*{Acknowledgements}
The author would like to thank the Reviewers and the Managing Editor of the journal \emph{Reports on Mathematical Physics} for their valuable comments and suggestions, which helped to improve the quality of the manuscript. Furthermore, the author would like to thank Professors Joseph Hoisington, Jun-ichi Inoguchi and Toru Sasahara  for helpful discussions on clarifying and correcting the previous version of this paper. Explicitely, compared with the previous version, we have harmonized the signs of two formulas taken from different references, but which corresponded to distinct sign conventions for the Riemannian curvature tensor. Consequently, the stability result was reformulated with the revised sign in mind.

Crina-Daniela NEAC\c{S}U\\ \emph{National University of Science and Technology Politehnica Bucharest,\\
Faculty of Applied Sciences,\\
Department of Mathematics and Informatics,\\
313 Splaiul Independen\c{t}ei, 060042 Bucharest, Romania}\\
E-mail: crina.neacsu@upb.ro
\end{document}